\documentclass[11pt]{article}

\usepackage{amssymb,amsfonts}

\newtheorem{thm}{Theorem}[section]
\newtheorem{lem}[thm]{Lemma}
\newtheorem{prop}[thm]{Proposition}
\newtheorem{cor}[thm]{Corollary}

 \newenvironment{pf}{\paragraph{Proof.}}{\par\medskip}

\newcommand{\qed}{\ifhmode\unskip\nobreak\fi\quad\ensuremath\diamond\medskip}

\newcommand{\C}{\mathbb C}
\newcommand{\pp}{\mathbb P}
\newcommand{\calO}{\mathcal{O}}

\title{ On surfaces with $p_g=2$, $q=1$ and non-birational bicanonical map }
 \author{Ciro Ciliberto, Margarida Mendes Lopes}
\date{}
\begin{document}

\maketitle
\begin{abstract}  In the  present
note we show that any  surface
of general type over $\C$ with $p_g=2$,$q=1$ and non birational
bicanonical map has a pencil of curves of
genus $2$. Combining this result with previous ones, one obtains that
an irregular surface $S$ of general type with $\chi(S)\geq
2$ and non-birational bicanonical map has a pencil of curves
of genus
$2$.

\end{abstract}
\section{Introduction}
If a smooth complex surface $S$ of general type has a pencil of curves of genus $2$, then the bicanonical
map $\phi$ of $S$ is not birational. We call this exception to the birationality of the bicanonical map
$\phi$ the {\it standard case}. The classification of the non-standard cases has a long history and we
refer to the expository paper \cite{ci} for information on this problem. \par

The classification of non-standard irregular surfaces has been first considered
by Xiao Gang in \cite{xiaocan}. He gives a list of
numerical possibilities for the invariants of the cases which might occur. More precise results have been
obtained in \cite {cfm}, \cite{ccm} and \cite {cm}. \par

In this note we prove the following:

\begin{thm} \label{mainthm} Let $S$ be an irregular surface of general type over $\C$ with $\chi(S)\geq 2$
and non-birational bicanonical map. Then $S$ has a pencil of curves of genus $2$.\end{thm}

In view of the results contained in the aforementioned papers, in order to prove this theorem, it suffices
to exclude the existence of non-standard cases with $p_g=2, q=1$. This we will do in \S \ref{theorem}
(see theorem \ref {main}). In \S \ref{preliminar} we prove some vanishing results
 which we apply later. In \S\S \ref {paracanonical}, \ref {nonbirat} we study the
paracanonical system of the surfaces in question, proving a few numerical and geometric properties
  of it which allow us to prove theorem \ref {main}.\bigskip

\noindent {\bf Acknowledgements:} The present collaboration took place in the framework of the EC
research project HPRN-CT-2000-00099, EAGER. \par

The second author is a member of CMAF and of the Departamento
de Matem\'atica da Faculdade de Ci\^encias da Universidade
de Lisboa. \par
The present paper has been finished during a visit of
the second author to Rome supported by GNSAGA of
INDAM.

Finally the authors want to thank the referee for the careful reading of the
paper and for useful comments.\bigskip

\section{A preliminary result}\label{preliminar}

In this section we prove some results that we will need  further on.

\begin {lem} \label{h}
Let $S$ be an irregular surface of general type.

 If $C$ is a $1$-connected curve on $S$ and $\eta\in Pic^0(S)$ is a point such that $h^1(S,{\cal
O}_S(K+\eta))=0$, then
$h^1(S,{\cal O}_S(K+\eta+C))\leq 1$ and equality holds if and only if $\eta_{|C}\simeq {\cal O}_C$.\par
Furthermore, if $h^1(S,{\cal O}_S(K+C))=0$ and $h^1(S,{\cal O}_S(K+\eta+C))=1 $, then $\eta$ is a torsion point of $Pic^0(S)$.
\end{lem}

\begin{pf}  By  Serre
duality and the hypothesis, we have $h^1(S,{\cal O}_S(K+\eta))=h^1(S, {\cal O}_S(-\eta))=0$.

Considering the long exact sequence obtained from the exact sequence:\par

$$0\to {\cal O}_S(-\eta-C)\to {\cal O}_S(-\eta)\to {\cal O}_C(-\eta)\to 0$$

\noindent  we obtain $h^1(S,{\cal O}_S(-\eta-C))= h^0(C,{\cal O}_C(-\eta))$. The line bundle ${\cal
O}_C(-\eta)$ has degree $0$ on
$C$, thus, by
\cite {cfm}, corollary (A.2),  $h^0(C,{\cal O}_C(-\eta))\leq 1$ with equality holding if and only if   $-\eta_{|C}\simeq
{\cal O}_C$. Now the first assertion
  follows because, by Serre duality $h^1(S,{\cal O}_S(-\eta-C))=h^1(S,{\cal O}_S(K+\eta+C))$.\par The second assertion
follows immediately from the first if we recall that
$h^1(S,{\cal O}_S(K+C))$ is the dimension   of the kernel of the restriction map $Pic^0(S)\to Pic^0(C)$.
\qed\end{pf}

\begin {prop} \label{h1}
Let $S$ be an irregular surface of general type such that either $q=1$ or the image of the
Albanese map $a:S\to Alb(S)$ is a surface. \par

\noindent (i) If $C$ is a $1$-connected curve on $S$ and $\eta\in Pic^0(S)$ is a general point, then
$h^1(S,{\cal O}_S(K+\eta+C))\leq 1$.\par

\noindent (ii) If $C$ is a $1$-connected curve on $S$ which is not contracted by the Albanese map $a:
S\to Alb(S)$, then $h^1(S,{\cal O}_S(K+C))<q$ and if $\eta\in Pic^0(S)$ is a general point, then
$h^1(S,{\cal O}_S(K+\eta+C))=0$.

\noindent (iii) If $C$ is a $1$-connected curve on $S$ which is contracted by the Albanese map $a:
S\to Alb(S)$, then $h^1(S,{\cal O}_S(K+C))=q$ and if $\eta\in Pic^0(S)$ is a general point, then $h^1(S,{\cal
O}_S(K+\eta+C))=1$.\end{prop}

\begin{pf}
\noindent  Assertion (i) is an immediate consequence of lemma \ref{h} because $h^1(S,{\cal O}_S(K+\eta))=0$ for a
general point
$\eta\in Pic^0(S)$. For $q=1$ this is clear by upper semicontinuity of $h^0(S,{\cal
O}_S(K+\eta))$ on $\eta\in Pic^0(S)$, whereas for $q\geq 2$ it follows
by the generic vanishing theorem from \cite {gl} (see also \cite {bea}).

Let us prove assertion (ii). The assertion about $h^1(S,{\cal O}_S(K+C))$ is well known (see \cite {ra} and \cite{ca},
remark 6.8). Indeed, since $C$ is not contracted by the Albanese map, if
$\eta\in Pic^0(S)$ is a general point, then $\eta_{|C}$ is non-trivial, i.e. the restriction map $Pic^0(S)\to Pic^0(C)$ is
non-zero.  Hence its tangent map $H^1(S,{\cal O}_S)\to H^1(C,{\cal O}_C)$ is non-zero, thus $h^1(S,{\cal O}_S(K+C))$, which
is the dimension of its kernel, is smaller than $h^1(S,{\cal O}_S)=q$.  The assertion about $h^1(S,{\cal O}_S(K+\eta+C))$ then
follows also from lemma \ref{h}.

The proof of (iii) is similar and therefore we omit it.\qed \end{pf}

\section{Some properties of the paracanonical curves}\label {paracanonical}

Let $S$ be an irregular surface of general type. We denote by $K$ a canonical divisor of $S$.
If $\eta\in Pic^0(S)$ is any point, we will consider the linear system $|K+\eta|$. We denote by $C_\eta$ a
curve in $|K+\eta|$, called a {\it paracanonical curve} of $S$. \par

Assume now that $p_g=2, q=1$ for $S$. Then,  for a general point $\eta\in
Pic^0(S)$, the linear system $|K+\eta|$ is a pencil. The curves $C_\eta$ thus describe, for a general
point $\eta\in Pic^0(S)$, a continuous system $\cal K$ of dimension $2$ of curves on $S$, called the
{\it main paracanonical system} of $S$.\par

We will write
$|K+\eta|=F_\eta+|M_\eta|$, where $F:=F_\eta$ is the fixed part and $|M|:=|M_\eta|$ is the movable
part of $|K+\eta|$. Next we prove two lemmas about the main paracanonical system.

\begin{lem} \label {vary} In the above setting, there is a not empty open Zariski subset $U\subseteq Pic^0(S)$ such
that either $F_\eta$ or $|M_\eta|$ stays fixed as $\eta$ varies in $U$. \end{lem}

\begin{pf} Let $U\subseteq Pic^0(S)$ be the not empty open set such that $F_\eta$ depends flatly on $\eta\in U$,
so that the algebraic equivalence class of $F_\eta$ is independent on $\eta\in U$.\par

Suppose $F_\eta$ varies with $\eta\in U$. Fix a point $\epsilon\in Pic^0(S)$. Then,
when $\eta$ varies, the system $F_\eta+|M_\epsilon|$ varies describing the general paracanonical system
in $\cal K$. This proves that $|M_\epsilon|$ does not depend on $\eta$. \qed \end {pf}

\begin{lem}\label{fm}  Let $S$ be a minimal surface of general
 type with $p_g=2$, $q=1$, with no pencil of curves of genus $2$ and $K^2\leq 8$. For any $\eta\in Pic^0(S)$ such that
$h^0(S, {\calO}_S(K+\eta))=2$, the general curve in  $|M_\eta|$ is irreducible. Furthermore:\par

\noindent (i) if $|M_\eta|=|M|$ does not depend on $\eta\in Pic^0(S)$, then either
$F_\eta$ is $1$-connected or $F_\eta=A_\eta+B$, where $B$ is a fundamental cycle, $A:=A_\eta$
is $1$-connected, and $K^2_S=8$, $K_S\cdot M=4$, $A\cdot B=0$, $K\cdot A=4$, $A^2=2$, $M\cdot
A=M\cdot B=2$;\par

\noindent (ii) if $F_\eta=F$ does not depend on $\eta\in Pic^0(S)$, then either $F=0$,
or $F\cdot M=2$, $F$ is $1$-connected and $M^2\geq 4$ or $K^2=8$, $M^2=F\cdot M=4$ and $K\cdot F=0$. In
this latter case either $F$ is $1$-connected or $F=A+B$ with $A$, $B$ fundamental cycles such that
$A\cdot B=0$.  \end{lem}

\begin{pf} We start by proving that the general curve in $|M_\eta|$ is irreducible if $h^0(S, {\calO}_S(K+\eta))=2$.
Otherwise the pencil $|M_\eta|$ would be composed with a pencil ${\cal P}$. The pencil ${\cal P}$ cannot be rational, since
$|M_\eta|$ has dimension $1$. Since $q=1$, then ${\cal P}$ would be an elliptic pencil, i.e. the Albanese
pencil of $S$, and $M_\eta\equiv 2G$, where $G$ is a curve in ${\cal P}$. Since $G^2=0$, we must have
$F_\eta\not=0$ and we cannot have $K\cdot F_\eta=0$, otherwise $F_\eta$ would be contained in curves of
${\cal P}$ and the canonical curves would be disconnected. Thus we have $8\geq K^2\geq K\cdot F+2K\cdot
G>2K\cdot G$, hence $K\cdot G=2$, i.e. ${\cal P}$ would be a pencil of curves of genus $2$, a
contradiction.\par

Now we prove (i). Suppose $|M_\eta|=|M|$ does not depend on $\eta\in Pic^0(S)$. Then $F_\eta$ moves with
$\eta$. We claim that $K\cdot F_\eta\geq 4$ and therefore $K\cdot M\leq 4$. Indeed, let $G:=G_\eta$ be an
irreducible component of $F:=F_\eta$ moving with $\eta$, so that $G$ is nef. Notice that by Debarre's inequality
\cite {debarre} we have $K^2\geq 2p_g=4$. Then the index theorem tells us
that $K\cdot G\geq 2$ and
$K\cdot G=2$ would imply that $G$ moves in a pencil of curves of genus $2$, a contradiction. If
$K\cdot G=3$, the index theorem again says that $G^2=1$ and we have a contradiction by
proposition (0.18) from \cite {ccm}. Therefore $K\cdot F\geq K\cdot G\geq 4$ and there is a unique
component $G$ of $F$ moving with $\eta$, otherwise $K\cdot F\geq 8$ and therefore $K\cdot M\leq 0$, a
contradiction.\par

If $F\cdot M=2$, then $F$ is $1$-connected by A.4 from \cite {cfm}. Assume $F\cdot M>2$. Then we
have $4\geq K\cdot M=F\cdot M+M^2\geq 4+M^2$, yielding $F\cdot M=4, M^2=F^2=0, K\cdot M=K\cdot F=4,
K^2=8$. If $F$ is not $1$-connected, then $F=A+B$ with $A, B$ effective, non-zero divisors, such that
$A\cdot B\leq 0$. Then $(M+A)\cdot B\geq 2$, hence $M\cdot B\geq 2-A\cdot B\geq 2$, and similarly
$M\cdot A\geq  2$. Since $F\cdot M=4$ we must have $M\cdot A=M\cdot B=2$ and $A\cdot B=0$. Thus $F$ is
$0$-connected and $A$ and $B$ are also $0$-connected by A.4 from \cite {cfm}. Suppose the unique
component $G$ of $F$ moving with $\eta$ sits in $A$. Then $K\cdot A=4$ and $K\cdot B=0$. From $0=K\cdot
B=B\cdot A+B^2+B\cdot M=B^2+2$, we deduce $B^2=-2$, hence $B$ is a fundamental cycle. Suppose $A$ is
not $1$-connected and write $A=A'+A''$, with $A'\cdot A''=0$. We claim that $B\cdot A'=B\cdot A''=0$.
Suppose indeed that $B\cdot A'>0$. Since $0=B\cdot  A=B\cdot A'+B\cdot A''$, we would have $B\cdot
A''<0$. But then $A''\cdot (A'+B)=A''\cdot B<0$, contrary to the fact that $F$ is $0$-connected.
Suppose $G$ sits in $A'$. We claim that $M\cdot G\geq 2$. Indeed, clearly $M\cdot G> 0$. If $M\cdot G=1$, since $G$ is
not rational, then $|M|$ would cut a fixed point on $G$, hence $G$ would be properly contained in a curve of $|M|$,
contrary to the fact that $G$ moves on $S$. Since $M\cdot G\geq 2$,
we have $M\cdot A'=2, M\cdot A''=0$. But then $A''\cdot (A'+B+M)=0$, contrary to the $2$-connectedness
of the canonical divisors. This ends the proof of (i).\par

In case (ii) we set
$M:=M_\eta$. One has $M^2>0$ since $M$ moves in a system ${\cal M}$ of curves of dimension $2$ and the general curve $M$
is irreducible. The case $M^2=1$ is excluded by proposition (0.14, iii) of \cite {ccm}. The case $M^2=2$ is also excluded
by theorem (0.20) of \cite {ccm}. Thus we may assume $M^2\geq 3$. Suppose that $M^2=3$, in which case $F\not=0$,
otherwise we would have $K^2=M^2=3$ against \cite {debarre}. Again by
theorem (0.20) of \cite {ccm}, ${\cal M}$ has no base point.  Fix
$\eta\in Pic^0(S)$ general and a general curve $M\in |M_\eta|$. For a general $\epsilon\in Pic^0(S)$,
the pencil $|M_\epsilon|$ cuts out on $M$ a base point free $g^1_3$ which has to vary with $\epsilon$ by
(1.6) of \cite {cfm}. This is impossible, since $M\cdot(K+M)= M\cdot (2M+F)\geq 8$ and
therefore $p_a(M)\geq 5$. Thus $M^2\geq 4$. If $M\cdot F\geq 4$, then $8\geq K^2\geq K\cdot M=F\cdot
M+M^2$, proving that $K^2=8$ and $M^2=F\cdot M=4$, $K\cdot F=0$. Suppose $F$ is not $1$-connected and
write $F=A+B$ with $A, B$ effective, non-zero divisors, such that $A\cdot B\leq 0$. By the same
argument we made before, we have $M\cdot A=M\cdot B=2$ and $A\cdot B=K\cdot A=K\cdot B=0$. Since $F^2=-4$ we have
that $A^2=B^2=-2$, and this concludes the proof of (ii).  \qed \end{pf}

\section{Non birationality of the bicanonical map and the paracanonical system}\label{nonbirat}

Now we will consider a minimal surface $S$ of general type with $p_g=2, q=1$, with non-birational
bicanonical map $\phi$, presenting the non-standard case.
Notice that, by
\cite {reider} and by proposition 3.1 of \cite {cm}, one has then
$K^2\leq 8$. We also recall the following lemma already contained in
\cite {ccm} (see lemma (2.2) of that paper):\par

\begin{lem}\label {lem}   Fix $\eta\in Pic^0(S)$ and
let $x$, $y$ be points on $S$ such that $\phi(x)=\phi(y)$, with $y$ not lying in the base
locus of $|K-{\eta}|$. In particular one may assume that $x$, $y$ are general
points such that $\phi(x)=\phi(y)$. Then $x$ belongs to a curve $C_\eta$ in $|K+\eta|$ if
and only if $y\in C_\eta$.  \end{lem}

Let $\nu>1$ be the degree of $\phi$. In view of this lemma, $\phi$ restricts to
the general curve $M_\eta$ to a map of degree $\nu$ to its image. In particular, if $\phi$ has degree $2$, then the
bicanonical involution fixes the curves $M_\eta$ and therefore also the curves $F_\eta$. \par

Now we are ready for the proof of the following proposition:\par

\begin{prop}\label {fmmore} Let $S$ be a minimal surface of general
 type with $p_g=2$, $q=1$, with non-birational
bicanonical map, presenting the non-standard case. Then, with the same notation as in lemma \ref {fm},
we have that $F=F_\eta$ is independent of $\eta\in Pic^0(S)$, and either:\par

\noindent (i) $F=0$ or $F$ is $1$-connected, strictly contained in a fibre of the Albanese pencil of $S$
and $F\cdot M=2$, $M^2\geq 4$, or;\par

\noindent (ii)  $F=A+B$ with $A$, $B$ fundamental cycles such that $A\cdot B=0$ and $M^2=F\cdot M=4$,
$K^2=8$.\par

In any case for a general $\eta\in Pic^0(S)$, the general curve $M\in |M_\eta|$ is
bielliptic and for a general $\epsilon\in Pic^0(S)$, the linear system $|M_\epsilon|$ cuts out on $M$ a
complete linear series whose movable part is a $g^1_4$ composed with the bielliptic involution on
$M$.\par

The bicanonical map has then degree $2$ onto its image, and the bicanonical involution $\iota$ acts on
the general curve $M$ as the bielliptic involution.\end{prop}

\begin{pf} First, assume that $F$ is non-zero, $1$-connected and not contained in a
fibre of the Albanese pencil. We claim that this case is not possible.\par

Notice that $h^1(S,{\cal
O}_S(2K-M_\eta))=h^1(S,{\cal
O}_S(K+F_\eta-\eta))=0$ for $\eta\in Pic^0(S)$ a general point. Indeed,
 if $|M_\eta|$ moves with $\eta$, then $F=F_\eta$
stays fixed by lemma \ref {vary}, and the assertion follows by proposition \ref{h1}, (ii). If
$|M_\eta|$ stays fixed while $\eta$ moves, the divisor
$F_\eta-\eta$  stays also fixed, and we obtain the assertion by  the second part of lemma
\ref{h}, because $h^1(S,{\cal
O}_S(-F_\eta))=0$ and $\eta$ is not a torsion point, for a general $\eta\in Pic^0(S)$ .\par

Thus if $\eta\in Pic^0(S)$ is general, the
bicanonical system $|2K|$ cuts out on the general curve $M\in |M_\eta|$ a complete, non special series
of degree $2g-2+F\cdot M$ and dimension $g-2+F\cdot M$, where we set $g:=p_a(M)$. By lemma \ref {lem}
this series is composite with an involution of degree $\delta\geq 2$. One has $2g-2+F\cdot M\geq \delta
(g-1)+\delta(F\cdot M-1)$, hence $F\cdot M\geq (\delta-2)(g-1)+\delta(F\cdot M-1)\geq 2(F\cdot M-1)$,
and therefore $F\cdot M=\delta=2$ and the equality has to hold everywhere in the above inequalities. In
particular $M$ is hyperelliptic. Notice that two points $x,y\in M$ are conjugated in the
hyperelliptic involution if and only if $\phi(x)=\phi(y)$, so that $\phi$ has degree $2$ onto its
image and we can talk about the bicanonical involution $\iota$ such that $\iota(x)=y$ if and only if
$\phi(x)=\phi(y)$.\par

Suppose that $|M_\eta|$ moves with $\eta$. Let $\epsilon\in Pic^0(S)$ be a general point and consider
the linear series $g_{\epsilon,\eta}$ cut out by $|M_\epsilon|$ on the general curve $M\in |M_\eta|$.
This series is complete and its movable part is composite with the hyperelliptic involution
on $M$ by lemma \ref {lem}. Furthermore its fixed divisor does not depend on $\epsilon$, since it is
supported at the finitely many points of $S$ which belong to every curve $M_\eta$. This would imply that
$g_{\epsilon,\eta}$ is independent on $\epsilon$, a contradiction to (1.6) of \cite {cfm}.\par

Suppose that $|M_\eta|$ does not move with $\eta$ and let again $M$ be the general curve in $|M_\eta|$.
Let $\epsilon\in Pic^0(S)$ be a general point. By lemma \ref {lem} the curve $F_\epsilon$ cuts out on $M$
a divisor fixed by the hyperelliptic involution. This would imply that the restriction map $Pic^0(S)\to Pic^0(M)$ is the
zero map, yielding $h^1(S, {\cal O}_S(-M))=1$, a contradiction. Our claim is thus proved. \par

Suppose next that $F$ is $1$-connected, contained in some fibre of the Albanese pencil
of $S$. Remark that $h^1(S, {\cal O}_S(-F))=1$ and also $h^1(S, {\cal O}_S(-F+\eta))=1$, for $\eta\in Pic^0(S)$ general
by proposition \ref {h1}. Now we will prove that (i) holds and therefore $|M_\eta|$ moves with $\eta\in Pic^0(S)$.\par

In this case $|2K|$ cuts out on $M$ an incomplete  series of degree $2g-2+F\cdot M$ and
dimension $g-3+F\cdot M$, which is composite with an involution of degree $\delta\geq 2$. We claim that
$\delta=2$. Suppose indeed that $\delta\geq 3$. Then $2g-2+F\cdot M\geq 3(g-3)+3F\cdot M$, i.e.
$g+2F\cdot M\leq 7$, which yields $g=3$, $F\cdot M=2$, $\delta=3$. Then $4=(K+M)\cdot M=F\cdot M
+2M^2=2+2M^2$, and we find $M^2=1$. By theorem (0.20) of \cite {ccm}, we see that $|M|$ cannot vary with
$\eta$, which means that $F$ has to. Namely $F$ has to be a fibre of the Albanese pencil. Then $K\cdot
F=F^2+F\cdot M=2$, and we find a contradiction. \par

 \ From $\delta=2$ we deduce $2(g-3+F\cdot M)\leq 2g-2+F\cdot M$, hence $F\cdot M\leq 4$.\par

 If $F\cdot
M=4$, then $M$ has to be hyperelliptic. Arguing as above, we see that $|M|$ cannot move with $\eta$,
hence $F$ has to, so that $F$ is a fibre of the Albanese pencil. In this case, by lemma \ref {lem}, the
bicanonical involution has to fix every curve of the Albanese pencil. Thus, if $x,y\in M$ are in the
$g_2^1$, i.e. if $\iota(x)=y$, then the fibre of the Albanese pencil through $x$ also contains $y$. This
would imply that the base of the Albanese pencil is rational, a contradiction.\par

 If $F\cdot
M=2$, then $K\cdot
F=F^2+F\cdot M\leq 2$ and therefore $F$ cannot move, otherwise the fibres of the Albanese pencil have
genus $2$. Then, by proposition \ref {fm}, $M^2\geq 4$. This ends the proof of (i). \par

Suppose now that $F$ is not $1$-connected. Then we will prove that $|M_\eta|$
 moves with $\eta\in Pic^0(S)$, thus (ii) will follow by lemma \ref{fm}, (ii).\par

If $|M|$ does not
depend on $\eta\in Pic^0(S)$, then one has a decomposition $F=A+B$ as in part (i) of lemma \ref {fm} and we claim that
$h^1(S, {\cal O}_S(-F))=1$. Indeed $h^1(S,{\cal O}_S(-A))=0$ and $h^0(B,{\cal O}_B(A))=1$ and the claim follows by the
exact sequence

$$0\to {\cal O}_S(-F)\to {\cal O}_S(-A)\to {\cal O}_B(-A))\to 0$$

\noindent Hence $|2K|$ cuts out on the general curve $M\in |M|$ a linear series of degree $2g+2$ and dimension
$r\geq g+1$ , which is composite with an involution of degree $\delta\geq 2$. This implies again that the involution in
question is a $g_2^1$ on $M$. Arguing as before we see that the curves $A$ should cut out on $M$ divisors of
this $g^1_2$, which leads to a contradiction. Thus we proved that $|M_\eta|$
 moves with $\eta\in Pic^0(S)$. \par

As for the final part of the statement, first suppose we are in case (i). As we saw, $|2K|$ cuts out on
$M$ a $g_{2g}^{g-1}$ composed with an involution of degree $2$. As usual, one sees that $M$ cannot be
hyperelliptic. Then it has to be bielliptic. The series cut out by $|M_\epsilon|$ on $M$ is complete.
Its fixed part does not depend on $\epsilon$, whereas its movable part is composed with the bielliptic
involution and it is complete, of dimension $1$, hence it is a $g_4^1$.\par

Suppose now we are in case (ii). Recall that, by lemma \ref {lem}, the degree of the map $\phi_M:M\to \phi(M)$ is $\nu\geq
2$, which is also the degree of $\phi$. Since $M\cdot (K+M)=12$, then $M$ has genus $7$. The series cut out by
$|M_\epsilon|$ on $M$ is a complete $g_4^1$, varying with the
parameter $\epsilon\in Pic^0(S)$. Notice that its fixed part does not depend on $\epsilon$, hence it is
empty, since there are finitely many $g^1_d$'s on $M$, with $d\leq 3$. Furthermore lemma \ref {lem} yields that either
$\nu=2$ or $\nu=4$. In the latter case however the $g_4^1$ would not vary with $\epsilon$, a contradiction. Therefore
$\nu=2$ and the $g_4^1$ is composed with a fixed involution of degree $2$, independent on $\epsilon$. Since
the $g_4^1$ is complete and varies with $\epsilon$, we see that the involution on $M$
has to be elliptic.\qed\end{pf}

Notice that the assertion concerning the degree of the bicanonical map contained in the above
proposition, also follows by the results of \cite {xiaocan}.\bigskip

\section{The main theorem}\label {theorem}

In this section we will prove the following:

\begin {thm}\label {main} If $S$ is a surface of general type with
$p_g=2$, $q=1$, non-birational bicanonical map, then it presents the standard case.\end {thm}

By proposition \ref {fmmore}, if $S$ is a minimal surface of general type with $p_g=2$, $q=1$,
such that the bicanonical map is not birational and $S$ presents the non-standard case, then for
$\eta\in Pic^0(S)$ a general point, the pencil $|K+\eta|=F+|M_\eta|$ has a fixed part $F$ which does not
depend on $\eta$ and a movable part $|M|=|M_\eta|$ which depends on $\eta$.  Notice that all components
of $F$ sit in fibres of the Albanese pencil.\par

We let ${\cal M}$ be the
$2$-dimensional system of curves which is the closure, in the Hilbert scheme, of the family of curves
$|M_\eta|$, with $\eta\in Pic^0(S)$ a general point. If we let $\mu:=M^2$, then $\mu\geq 4$ and
the general curve $M\in {\cal M}$ is irreducible. Notice that ${\cal M}$ has a natural morphism
${\cal M}\to Pic^0(S)$ whose general fibre is a $\pp ^1$. \par

Fix a general point $x\in S$ and set
$x':=\iota(x)$. Consider the system ${\cal M}_{x}$ of curves in ${\cal M}$ passing through $x$. We claim
that the system ${\cal M}_{x}$ is irreducible, of dimension $1$, parametrized by
$Pic^0(S)$. Indeed any irreducible component of ${\cal M}_{x}$ has dimension $1$ and no irreducible
component of ${\cal M}_{x}$ can be contained in a fibre of ${\cal M}\to Pic^0(S)$. Since ${\cal M}_{x}$
cuts the general fibre of ${\cal M}\to Pic^0(S)$ in one point, the claim follows.\par

Notice that ${\cal M}_{x}={\cal M}_{x'}$, so it is appropriate to denote this
system by ${\cal M}_{x,x'}$.\par

 Let us point out the following corollary of
proposition \ref {fmmore}:\par

\begin {cor} \label {corol} In the above setting, two general curves in ${\cal M}$ have intersection
multiplicity $\mu-4$ at fixed points of $S$, and intersect at $4$ further distinct variable points,
which are pairwise conjugated in the bicanonical involution $\iota$. In particular, given a general
point $x\in S$, two general curves in ${\cal M}_{x,x'}$ intersect at $2$ distinct variable
points, which are conjugated in the bicanonical involution.\end {cor}

Recall that the {\it index} of a $1$-dimensional system of curves on
a surface is the number of curves of the system passing through a general point of the surface. Next we
prove the following lemma:

\begin{lem}\label {index} In the above setting, the system ${\cal M}_{x,x'}$ has index
$\nu=2$.\end{lem}

\begin{pf}  The index of ${\cal M}_{x,x'}$ cannot be $1$. In this case, in fact, ${\cal
M}_{x,x'}$ would be a pencil. Since ${\cal M}_{x,x'}$ has base points at $x$ and $x'$, it would be a
rational pencil, whereas we know it is parametrized by $Pic^0(S)$. Thus $\nu\geq 2$. \par

Let $y$ be a general point of $S$ and set $y':=\iota(y)$. Let $M$ be a curve in ${\cal M}_{x,x'}$
through $y$ and therefore also through $y'$. Since $\nu\geq 2$, we know there is some other curve $M'$
in ${\cal M}_{x,x'}$ through $y$ and  $y'$. Suppose there is a third one $M''$. Then $M'$ and $M''$ would
cut out on $M$ the same divisor, a contradiction to (1.6) of \cite {cfm}. \qed \end{pf}

Set $Pic^0(S):=A'$ and fix $x\in S$ a general point. Another general point $y\in S$ determines two point
$m_{1,y}, m_{2,y}$ in $A'$ corresponding to the two curves of ${\cal M}_{x,x'}$ containing
$y$. Thus we can consider the map:\par

$$\alpha: S\to Pic^2(A')\simeq A'$$

\noindent which takes the general point $y\in S$ to the divisor class of $m_{1,y}+m_{2,y}$. This map is
clearly surjective and it factors through the Albanese map. We denote by $G$ the general fibre of
$\alpha$, which is composed of curves of the Albanese pencil. \par

We are finally in a position to give the:\medskip

 \noindent {\bf Proof of theorem \ref{main}}. Assume $S$ presents the non-standard case. Let us keep the
above notation and let us set $n:=G\cdot M$.  By proposition \ref {fmmore}, one has $G\cdot K=G\cdot
F+G\cdot M=G\cdot M=n$. Since we are in the non-standard case, we have $n\geq 4$.\par

Fix a general point $x\in X$ and the system
${\cal M}_{x,x'}$. Let $M$ be a general curve in ${\cal M}_{x,x'}$. Let $M\cap G$ consist of the points
$x_1,...,x_n$. Notice that each of the points $x_1,...,x_n$ is a general point of $S$, in particular it
is different from $x$ and $x'$. Let $m$ be the point of $A'$ corresponding to the curve $M$. By lemma
\ref {index} and by the generality of
 $x_1,...,x_n$, for each $i=1,...,n$ there is only another curve $M_i\in {\cal M}_{x,x'}$,
different from $M$, containing $x_i$. Let $m_i$ be the point of $A'$ corresponding to $M_i$,
$i=1,...,n$. One has $\alpha(x_i)=m+m_i$, $i=1,...,n$. On the other hand, by the meaning of $G$,
one has $\alpha(x_1)=...=\alpha(x_n)$, i.e. the divisor classes of $m+m_i$ on $A'$ are the same for all
$i=1,...,n$. This implies $m_1=...=m_n$, i.e. $M_1=...=M_n$. Let $M'$ be this curve. Then $x_1,...,x_n$
sit in the intersection of the two curves $M$ and $M'$ of ${\cal M}_{x,x'}$, off the points $x, x'$. By
corollary \ref {corol} we have $n+2\leq 4$, i.e. $n\leq 2$, a contradiction. \qed

\bigskip

\bigskip

\begin{tabbing}  Universit\`a di Roma Tor Vergata xxxxxxxxxxxxxxxx\=Universidade de Lisboa  \kill Ciro
Ciliberto \> Margarida Mendes Lopes\\
 Dipartimento di Matematica \> CMAF\\
 Universit\`a di Roma Tor Vergata \>Universidade de Lisboa \\ Via della Ricerca Scientifica \>
Av. Prof. Gama Pinto, 2 \\ 00133 Roma, ITALY \> 1649-003 Lisboa, PORTUGAL\\
cilibert@mat.uniroma2.it \> mmlopes@lmc.fc.ul.pt
\end{tabbing}

\end{document}